\documentclass[11pt]{article}

\title{Vassiliev theory and regional change}

\author{Jonathan Fine\relax
\thanks{203 Coldhams Lane, Cambridge, CB1 3HY, England.  
\quad E-mail: \texttt{j.fine@pmms.cam.ac.uk}\hfill\break
(First version: 26 November 1996)
Postscript (2 March 1998).  An expanded version of this paper is in
preparation.}
}
\date{7 December 1996}

\def\bfR{{\bf R}}
\def\Kcal{{\cal K}}
\newcommand\set[1]{^{\{#1\}}}
\newcommand\ord[1]{^{(#1)}}
\newcommand\Kquot[2]{\Kcal #1{#2}/\Kcal #1{#2+1}}
\newcommand\Kbar{\smash{\overline{K}}\vphantom{K}}

\begin{document}
\maketitle

The purpose of this note is to state some definitions that may be useful
in the study of knots, manifolds and the like.  They apply to anything for
which the concept of a regional change can be defined, such as a product
of elements in a group.

The motivation comes from the Vassiliev theory for invariants of knots in
$\bfR^3$, the general nature of the axioms for a topological quantum field
theory, and an observation regarding the Kontsevich knot integral.  The
basic concept is that of a regional change. 

To change a crossing in a knot diagram is an example of a regional change.
This provides a choice.  Assume that one option is labelled $a$ (for
above) and the other $b$ (for below).  Call such a double point.  (This is
not quite the standard usage.) Let $k$ be a knot with $r$ double points. 
The double points are assumed ordered, and labelled $1$ to $r$.  Now for
each word $w$ in $a$ and $b$ of length $r$ let $k^w$ be the actual knot
formed by taking the $w_i$ choice at the $i$-th double point.  

Let $\Kcal$ denote the group of formal sums of (isomorphism types of)
knots in $\bfR^3$, with integer coefficients.  Each knot $k$ with $r$
regional changes (or double points) determines a formal sum
\[
    \{k\} = \sum \nolimits _{|w|=r} (-1)^w k^w
\]
of knots.  Here $(-1)^w$ is $(-1)^n$, where $n$ is the number of $b$'s in
$w$.  The formal sum
\[
    (k) = \sum \nolimits _{|w|=r} w k^w
\]
of $w$-weighted knots is also useful.  Note that $\{k\}$ does not depend
on the order of the double points, and only up to a sign on the $a$ and
$b$ labelling.  Define $\Kcal\set r$ and $\Kcal \ord r$ to be the groups
generated by all such $\{k\}$ and $(k)$ respectively.  The first is a
subgroup of $\Kcal$, the second of $\Kcal^N$, where $N=2^r$.

A Vassiliev invariant of order $r$ is simply a homomorphism from the
quotient group $\Kbar\set r = \Kquot \set r$.  It is easy to show that
these groups are finitely generated, but determination of their rank is
difficult.  They have been much studied.

The groups $\Kcal\ord r$ can be assembled into a complex.  Define
$\partial_i w$ to be the result of omitting the $i$-th letter from $w$,
multiplied by $(-1)$ if that letter is $b$.  Treat knots as constants.
The alternating sum
\[
    d = \partial_1 - \partial_2 + \partial_3 - \ldots 
        + (-1)^{r+1}\partial_r
\]
is the boundary map.  The resulting homology $\Kbar\ord r =
\Kcal\ord r / d\Kcal\ord{r+1}$ will be called the $r$-th difference group
(for knots).  Probably, these groups are finitely generated. (There
appears to be an argument involving permutations that establishes this, at
no extra cost over the Vassiliev proof.)  They seem to give more
information.

Now let $m$ be a manifold, and suppose that on $m$ there are given $r$
submanifolds with boundary.  Assume they are disjoint.  Call them $m_i^a$.
These are the regions.  Now for each $i$ provide $m_i^b$, which is to be a
manifold with boundary.  Assume that the $a$ and $b$ manifolds are
identified around their boundary.  It then follows that for each $w$ there
is a manifold $m^w$. Call such a system a manifold with $r$ regional
changes.  The quantities $\{m\}$ and $(m)$ can be defined just as before. 
This allows one to talk about Vassiliev invariants for manifolds, and also
difference groups.

Now let $k$ be a knot in a 3-manifold.  There is a fundamental group, so
crossing changes will not be enough to produce finitely generated
quotients.  Instead, the manifold definition should be imitated.  This
will allow the knot to be changed by an element in the fundamental group.

The Kontsevich integral expresses the order $r$ Vassiliev invariants as an
integral over $r$ horizontal slices moving up and down a knot $k$ in
$\bfR^3$.  This fact forces the following.  If within each of $r+1$ fixed
slices one makes a regional change to $k$, of any nature whatsoever, the
order $r$ Kontsevich integral vanishes on the alternating sum of the knots
so obtained.

The author has shown that the Vassiliev theories for 3-manifolds, and for
knots in 3-manifolds, are finitely generated.  There is reason to hope
that the same is true for all manifolds, and that the corresponding
invariants are of interest.

\end{document}